\documentclass[10pt,a4paper,oneside,leqno]{article}
\usepackage{amsfonts,amssymb}
\usepackage{eufrak}
\usepackage{amscd}

\newtheorem{defn}{Definition}

\vspace{0.5cm}

 4
\font\ebf=cmbx8
\font\erm=cmr8

\usepackage{graphicx, graphics}

\begin{document}

\thispagestyle{empty}

\begin{center}
	\noindent { \textsc{ On cobweb posets and their combinatorially admissible sequences}}  \\ 
	\vspace{0.3cm}
	\vspace{0.3cm}
	\noindent Andrzej Krzysztof Kwa\'sniewski \\
	\vspace{0.2cm}
	\noindent {\erm Member of the Institute of Combinatorics and its Applications  }\\
{\erm High School of Mathematics and Applied Informatics} \\
	{\erm  Kamienna 17, PL-15-021 Bia\l ystok, Poland }\noindent\\
	\noindent {\erm e-mail: kwandr@gmail.com}\\
	\vspace{0.4cm}
\end{center}

\vspace{0.4cm}

\noindent {\ebf SUMMARY}

\vspace{0.1cm}

\noindent {\small The purpose of this article is to pose three  computational
problems which are quite easily formulated for the new class of directed acyclic graphs
interpreted as Hasse diagrams. The problems posed are not yet solved though are of crucial 
importance  for the vast class of  new partially ordered sets with joint combinatorial interpretation.
These so called cobweb posets - are relatives of Fibonacci tree and are labeled by  specific
number sequences - natural numbers sequence and Fibonacci sequence included.  One presents here also a  
join combinatorial interpretation of those posets` $F$-nomial coefficients which are computed
with the so called cobweb admissible sequences. Cobweb posets and their natural subposets are 
graded posets, sometimes called a ranked posets. They are vertex partitioned into such antichains $\Phi_n$ 
(where $n$ is a nonnegative integer) that for each $\Phi_n$, all of the elements
covering $x$ are in $\Phi_{n+1}$ and all the elements covered by $x$ are in $\Phi_n$. We shall call
the $\Phi_n$ the  $n-th$- level. The cobweb posets may be identified with a chain of di-bicliques i.e. 
by definition - a chain of complete bipartite one direction digraphs.
Any chain of relations is therefore obtainable from the cobweb poset chain of complete relations
via deleting arcs in di-bicliques of the complete relations chain.}

\vspace{0.3cm}

\noindent Key Words: acyclic digraphs, tilings, special number sequences, binomial-like coefficients.

\vspace{0.1cm}

\noindent AMS Classification Numbers: 06A07 ,05C70,05C75, 11B39.

\vspace{0.2cm}

\noindent  affiliated to The Internet Gian-Carlo Polish Seminar:

\noindent \emph{http://ii.uwb.edu.pl/akk/sem/sem\_rota.htm} 

\vspace{0.1cm}

\noindent \textbf{Published} in : Adv. Studies Contemp. Math.  Vol.  \textbf{18}  No 1,  (\textbf{2009}),  17-32.

\vspace{0.4cm}

\section{Introduction}
In computer science a directed acyclic graph, also called DAG,
is a directed graph with no directed cycles. Utility of DAG`s is 
widely known and appreciated. Some algorithms - for example,
search algorithms - become simpler when applied to DAGs. 
DAGs  considered  as a generalization of trees  have a lot of applications in computer science,
bioinformatics,  physics and many natural activities of humanity and nature. 
For example in information categorization systems, such as folders in a computer or in
Serializability Theory of Transaction Processing Systems and many others.
Here we introduce specific DAGs as generalization of trees being
inspired by algorithm of the Fibonacci tree growth. For any given natural numbers valued sequence
the graded (layered) cobweb posets` DAGs  are equivalently representations of a chain of binary relations. Every relation of the cobweb poset chain is biunivocally represented by the uniquely designated  \textbf{complete} bipartite digraph-a digraph which is a di-biclique  designated by the very  given sequence. The cobweb poset is then to be identified with a chain of di-bicliques i.e. by definition - a chain of complete bipartite one direction digraphs. Any chain of relations is therefore obtainable from the cobweb poset chain of complete relations  via deleting  arcs (arrows) in di-bicliques.\\ 
Let us underline it again : \textit{any chain of relations is obtainable from the cobweb poset chain of  complete relations via deleting arcs in di-bicliques of the complete relations chain.} For that to see note that any relation  $R_k$ as a subset of  $A_k \times  A_{k+1}$ is represented by a
one-direction bipartite digraph  $D_k$.  A "complete relation"  $C_k$ by definition is identified with its one direction di-biclique graph  $d-B_k$
Any  $R_k$ is a subset of  $C_k$. Correspondingly one direction digraph  $D_k$ is a subgraph of an one direction digraph of  $d-B_k$.\\
The one direction digraph of  $d-B_k$ is called since now on \textbf{the di-biclique }i.e. by definition - a complete bipartite one direction digraph.
Another words: cobweb poset defining di-bicliques are links of a complete relations' chain. 

\vspace{0.2cm}

\noindent Because of that the cobweb posets in the family of all chains of relations unavoidably is of principle 
importance being the most overwhelming case of relations' infinite chains  or their finite parts (i.e. vide - subposets). 

\noindent The intuitively transparent names used above (a chain of di-cliques etc.) are to be the names of  defined 
objects in what follows. These are natural correspondents to their undirected graphs relatives
as we can view a directed graph as an undirected graph with arrowheads added. 

\noindent The purpose of this note is to put several questions intriguing on their own
apart from being fundamental for a new class of DAGs introduced below. 
Specifically this concerns problems which arise naturally in connection with
a new join combinatorial interpretation of all classical $F-nomial$
coefficients -  Newton binomial, Gaussian $q$-binomial and Fibonomial
coefficients included. This note is based on  [8,9] from which definitions and description
of these new DAG's are quoted for the sake of self consistency. 
\noindent Applications of new cobweb posets` originated Whitney numbers from [9] 
such as extended Stirling or Bell numbers are expected to be of at least such
a significance in applications to linear algebra of formal series 
as Stirling numbers, Bell numbers or their $q$-extended correspondent already are
in the so called  coherent states physics (see [13] for abundant references on the subject). 
Quantum coherent states physics is of course a linear theory with its principle of states` superposition.

\vspace{0.2cm}

\noindent \textbf{The problem to be the next.}
\noindent As cobweb subposets $P_n$ are   vertex partitioned into  antichains $\Phi_r$ for $r = 0, 1, ...,n$
which we call levels -  a question of canonical importance arrises.
\noindent Let $\left\{P_n\right\}_{n\geq0}$  be the sequence of finite cobweb subposets (see- below). 
What is the form and properties  of $\left\{P_n\right\}_{n\geq0}$'s  characteristic polynomials $\left\{\rho_n(\lambda)\right\}_{n\geq0}$  [15,17]?
\noindent For example - are these related to umbra polynomials? What are recurrence relations defining the $\left\{\rho_n(\lambda)\right\}_{n\geq0}$ family ?

\vspace{0.2cm}

\noindent \textbf{1.1.} Partially Ordered Sets - Elementary Information.

\noindent Let us recall indispensable notions thus establishing
notation and terminology.

\vspace{0.1cm}

\begin{defn}.
A partially ordered set (poset) is  an ordered pair $\langle P,
\leq \rangle$, where $P$ is a set and $\leq$ is a partial order on
$P$. Naturally $x<y  \Leftrightarrow  x\leq y  \wedge  x \neq y$.
\end{defn}

\begin{defn}.
An element $z$  covers an element $x$ provided that there exists
no third element  $y$ in between i.e. such $y$ that $x < y < z$.
In this case $z$ is called an upper cover of $x$.
\end{defn}
\noindent A graphical rendering of a partially ordered set (poset)
is being displayed via the cover relation of the partially ordered
set (poset) with an implied upward orientation by convention. It represents
the Hasse Diagram of  the poset [3].

\vspace{0.1cm}

\begin{defn}.
Hasse diagram is a graph constructed as follows. A point is drawn
for each element of the set $P$  and line segments are drawn
between these points according to the two rules: \textbf{1.} If $x
\leq  y $ then the point corresponding to $x$   appears lower in
the drawing than the point corresponding to $y$. \textbf{2.} The
edge between the points corresponding to elements $x$ and $y$ is
included in the graph if and only if   $x$  covers $y$ or $y$
covers $x$.
\end{defn}

\vspace{0.1cm}

\noindent In order to make the reading interactive one recommends
to go through preliminary guiding exercises  (cf. [3,2]).
\vspace{0.1cm}

\noindent \textbf{Exercise 1.} Draw  the Hasse diagrams for
Boolean algebras of orders n =2, 3, 4, and 5.

\vspace{0.1cm}

\noindent \textbf{Exercise 2.} Invent Hasse diagram for Fibonacci poset
$\langle P, \leq \rangle$. Its subsequent generations
$\Phi_s$ of elements $x,y\in P$ are designated by Fibonacci sequence as 
$x<y$ iff if $x$ is an ancestor of $y$. $F_0 = 0$ corresponds to the empty root $\{\emptyset\}$.

\vspace{0.1cm}

\noindent \textbf{Example}. Gallery of Posets by Curtis Greene may
be found in [2]. Compare therein the bottom levels of a Young-Fibonacci lattice, 
introduced by Richard Stanley with cobweb posets defined here below.

\vspace{0.1cm}

\noindent \textbf{1.2. Computation and Characterizing Problems}.
 
\noindent In the next section we define cobweb posets and their 
examples  are given [8,9]. A join combinatorial interpretation of cobweb
posets` characteristic binomial-like coefficients is provided too.
This simultaneously means join combinatorial interpretation of fibonomial coefficients and 
all incidence coefficients of reduced incidence algebras of full binomial type [16]. 

\noindent Finally we formulate three problems (to be explained
on the way):  characterization and/or computation of cobweb admissible sequences \textbf{Problem 1}, 
cobweb layers partition characterization and/or computation \textbf{Problem 2} and the
GCD-morphic sequences characterizations and/or computation \textbf{Problem 3} 
- all three  interesting on their own.

\section{Cobweb posets - presentation, examples and combinatorial interpretation}

\vspace{1mm}

\textbf{2.1. $F-nomial$ coefficients}.

\noindent Given any sequence $\{F_n\}_{n\geq 0}$ of nonzero reals
($F_0 = 0$ being sometimes acceptable as $0! = F_0! = 1.$)
one defines  its corresponding  binomial-like $F-nomial$
coefficients as in Ward`s Calculus of sequences [18] as follows.

\begin{defn}.$(n_{F}\equiv F_{n}\neq 0,\quad n > 0)$
$$
\left( \begin{array}{c} n\\k\end{array}
\right)_{F}=\frac{F_{n}!}{F_{k}!F_{n-k}!}\equiv
\frac{n_{F}^{\underline{k}}}{k_{F}!},\quad \quad n_{F}!\equiv n_{F}(n-1)_{F}(n-2)_{F}(n-3)_{F}\ldots 2_{F}1_{F};$$
$$ 0_{F}!=1;\quad n_{F}^{\underline{k}}=n_{F}(n-1)_{F}\ldots (n-k+1)_{F}. $$
\end{defn}
\noindent We have made above  an analogy driven identifications in the spirit of  Ward`s Calculus of sequences [18]. 
Identification $n_{F}\equiv F_{n}$ is  the notation used in extended Fibonomial Calculus case [10,11,12,6] being also there inspiring as $n_F$ mimics $n_q$ established notation for Gassian integers exploited in much elaborated family of various applications including quantum physics (see [10,9,13] and references therein).

\vspace{2mm}

\noindent \textbf{2.2. Cobweb infinite posets}. 
\noindent Cobweb infinite posets  $\Pi$  are designated uniquely by
any sequence of natural numbers $F=\{n_F\}_{n\geq 0}$ with the one root convention ($F_0=1$)  and are by construction endowed with a kind of  self-similarity property.
(In Fibonacci case  $F_0 = 0$ corresponds to the empty root $\{\emptyset\}$).
\noindent Given any such sequence $\{F_n\}_{n\geq 0}\equiv\{n_F\}_{n\geq 0}$ of 
positive integers  we define the partially ordered, graded infinite poset $\Pi$ - called a cobweb poset - as follows. Its vertices are labeled  by pairs of coordinates: ${\langle i , j \rangle} \in {N \times N_0}$ where $N$ stays for natural numbers while $N_0$ denotes the nonnegative integers.  Vertices show up in layers  of $N \times N_0$ grid along the recurrently emerging subsequent $s-th$ levels $\Phi_s$  ("generations") where  $s\in N_0$. 

\begin{defn}.
$
\Phi_s =\{\langle j, s\rangle 1\leq j \leq s_F\}, {s\in N_0}.
$
\end{defn}

\begin{defn}.
$$\Pi=\langle P,E\rangle,\  P=\bigcup_{0\leq p}\Phi_p ,\  E
=\{\langle\langle j , p\rangle ,\langle q ,(p+1) \rangle
\rangle\}\bigcup\{\langle\langle 1 , 0\rangle ,\langle 1 ,1
\rangle \rangle\},$$ 
$1 \leq j \leq {p_F} , 1\leq q \leq {(p+1)_F}.$
\end{defn}
We shall refer to $\Phi_s$  as to  the set of vertices at the
$s-th$ level. The population of the  $k-th$ level ("generation" )
counts  $k_F$  different member vertices for $k>0$ and one for
$k=0$. \vspace{2mm} Here down a disposal of vertices on $\Phi_k$
levels is visualized for the case of Fibonacci sequence.
$F_0 = 0$ corresponds to the empty root $\{\emptyset\}$.
On the account of the above definition we shall sometimes irrelevantly represent
$\Pi=\langle P,E\rangle,$ just by its vertex set  $P$ partitioned into countable family
of antichains (levels) uniquely designated by the $F$-sequence choice.  The same facility
is to be used for subposets on the grounds of a legal abuse of perfectionism.

\vspace{2mm}

$---- and ----- so ---- on ---- up    --- \Uparrow ----------$\\
$\star \star \star \star \star \star \star \star \star \star \star \star \star \star \star \star \star \star\star \star \star \star \star \star \star \star \star \star \star \star \star \star \star \star \star \star \star \star \star \star \star \star \star \star --\star \star \star \star\star10-th-level$\\
$\star \star \star \star \star \star \star \star \star \star \star \star \star \star \star \star \star \star \star \star \star \star \star \star \star \star \star \star \star \star \star\star\star\star----------- 9-th-level$\\
$\star \star \star \star \star \star \star \star \star \star \star \star \star \star \star \star \star \star \star \star \star------------------8-th-level$\\
$\star \star \star \star \star \star \star \star \star \star \star \star \star -------------------------7-th-level$\\
$\star \star \star \star \star \star \star \star-----------------------------6-th-level$\\
$\star \star \star \star \star ---------------------------------5-th-level$\\
$\star \star \star ---------------------------------- 4-th-level$\\
$\star \star -----------------------------------3-rd-level  $ \\
$\star ------------------------------------ 2-nd-level$\\
$\star ----------------------------------- 1-st-level$\\
$\emptyset ------------------------------------ 0-th-level$\\

\textbf{Figure 0. The $s-th$ levels in $ N\times N_0 $}

\vspace{1mm}

\noindent Accompanying to the set $P$ of vertices the set $E$ of edges
- we obtain the Hasse diagram. We may picture the partially ordered infinite set
$P$ from the  Definition $6$ with help of the sub-poset set $P_{m}$
({\it rooted at $F_{0}$ level of the poset}) to be continued then
ad infinitum in now obvious way as seen from the figures  $Fig.1-
Fig. 5$ of $P_{m}$ cobweb posets below. These non-trees look like 
Fibonacci-type trees \textit{with} a specific ``cobweb''.

\begin{center}

\includegraphics[width=75mm]{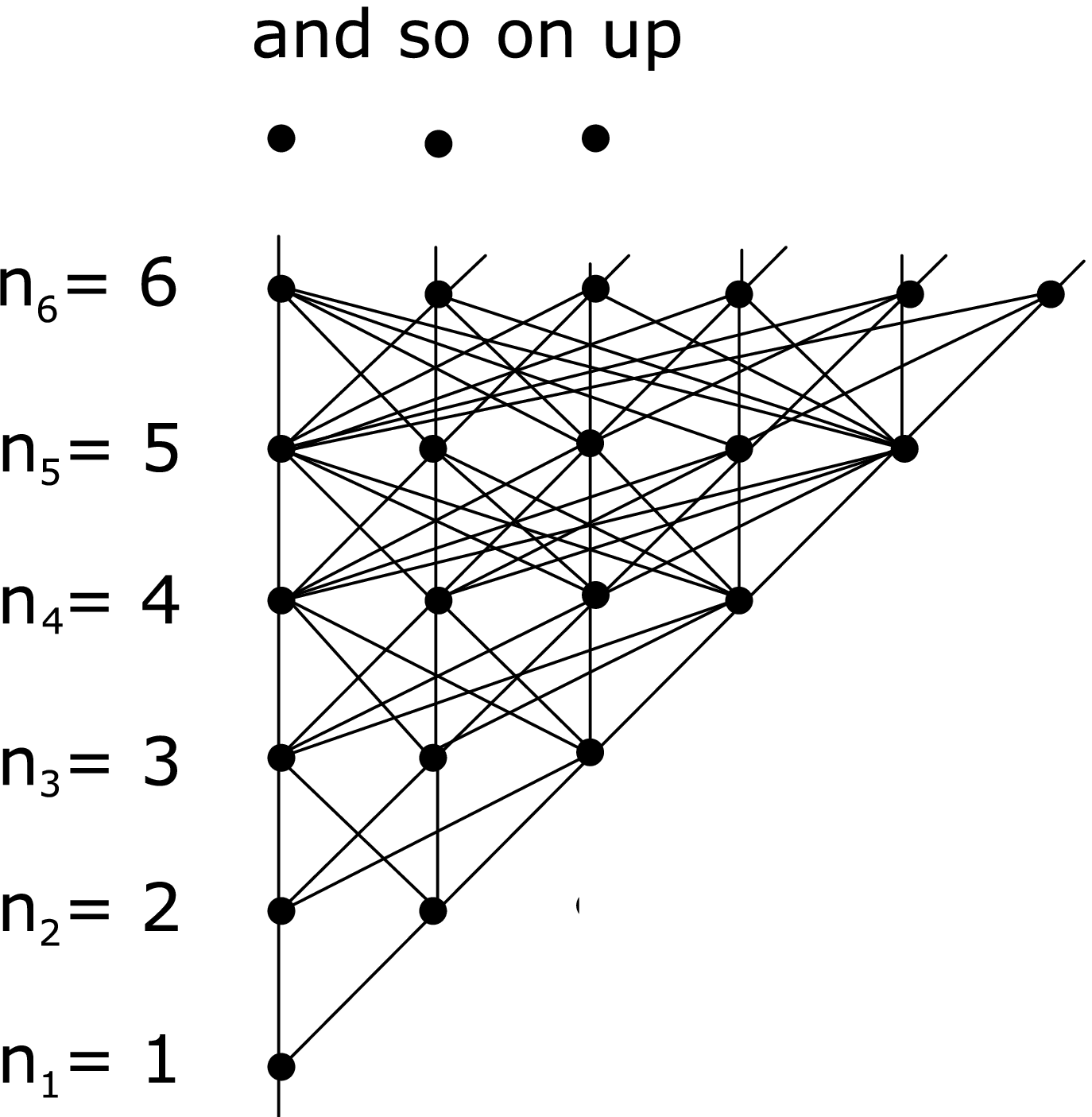}

\vspace{1mm}

\noindent {\small Fig.1. Display of Natural numbers' cobweb
poset.}

\end{center}

\vspace{1mm}

\begin{center}
\includegraphics[width=75mm]{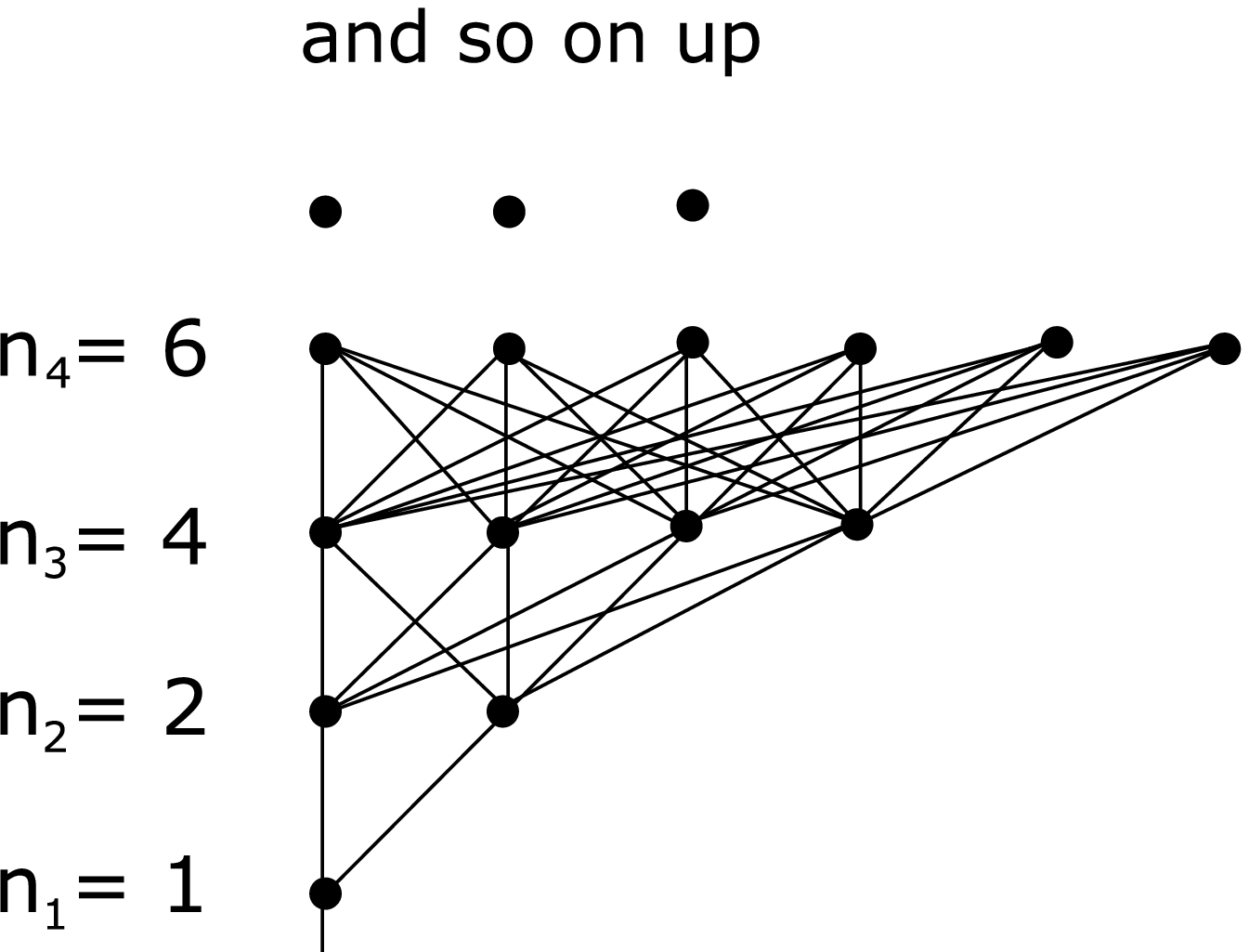}

\noindent {\small Fig.2. Display of Even Natural numbers' cobweb
poset.}
\end{center}

\vspace{1mm}

\begin{center}

\includegraphics[width=75mm]{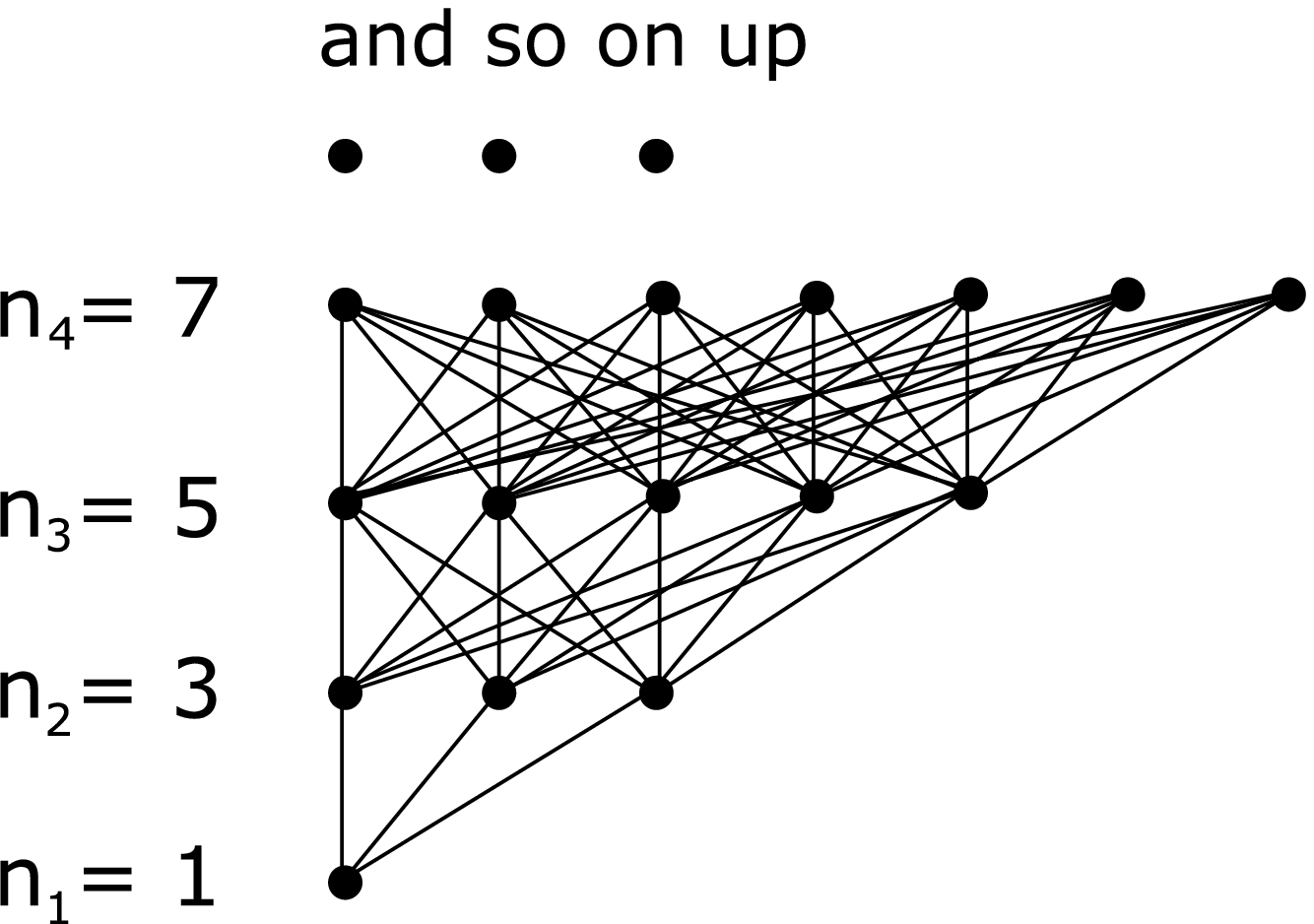}

\vspace{1mm}

\noindent {\small Fig3. Display of Odd natural numbers' cobweb poset.}
\end{center}

\vspace{1mm}

\begin{center}

\includegraphics[width=75mm]{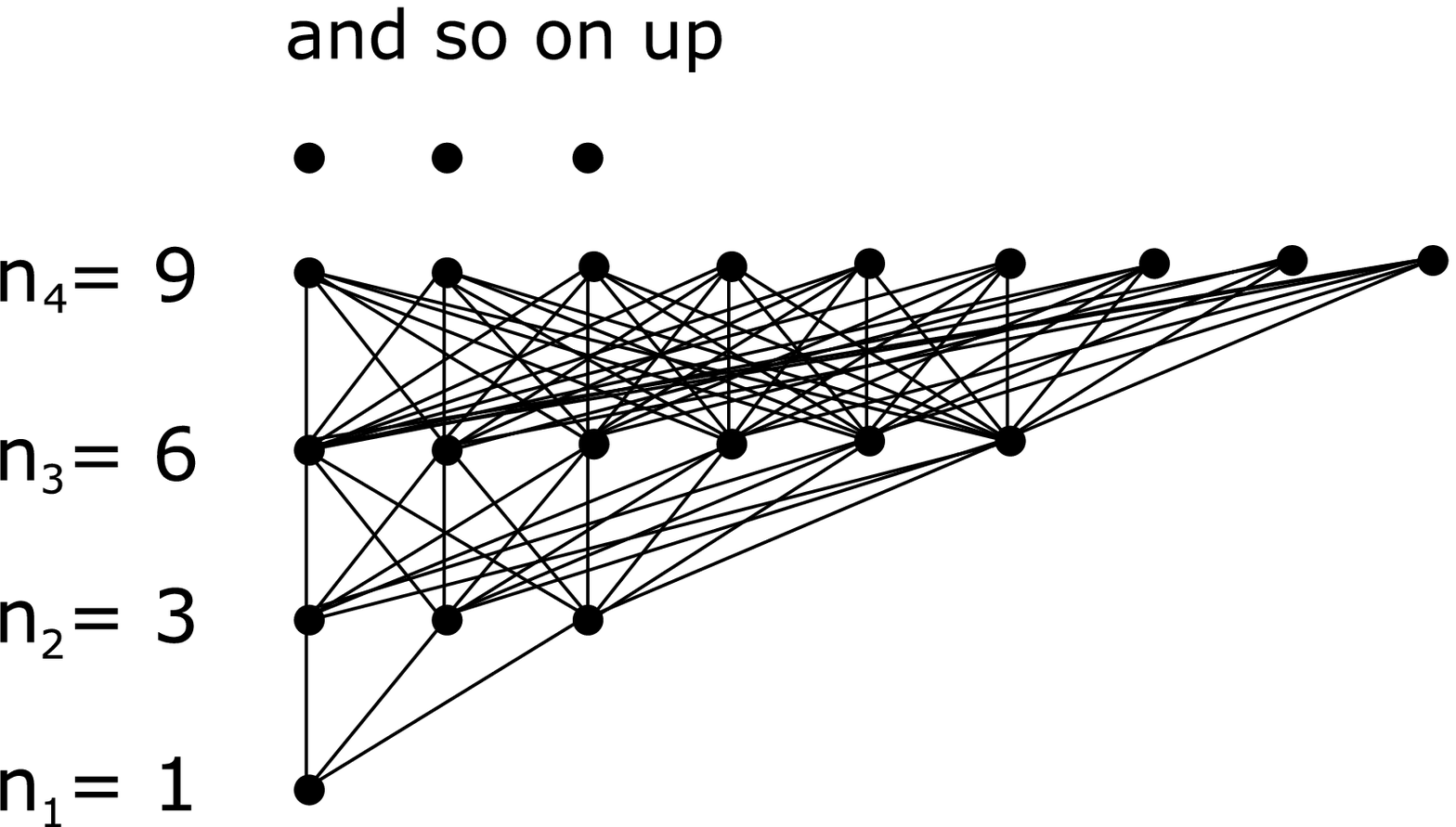}

\vspace{1mm}

\noindent {\small Fig.4. Display of divisible by 3 natural numbers' cobweb poset.}
\end{center}

\vspace{1mm}

\begin{center}

\includegraphics[width=75mm]{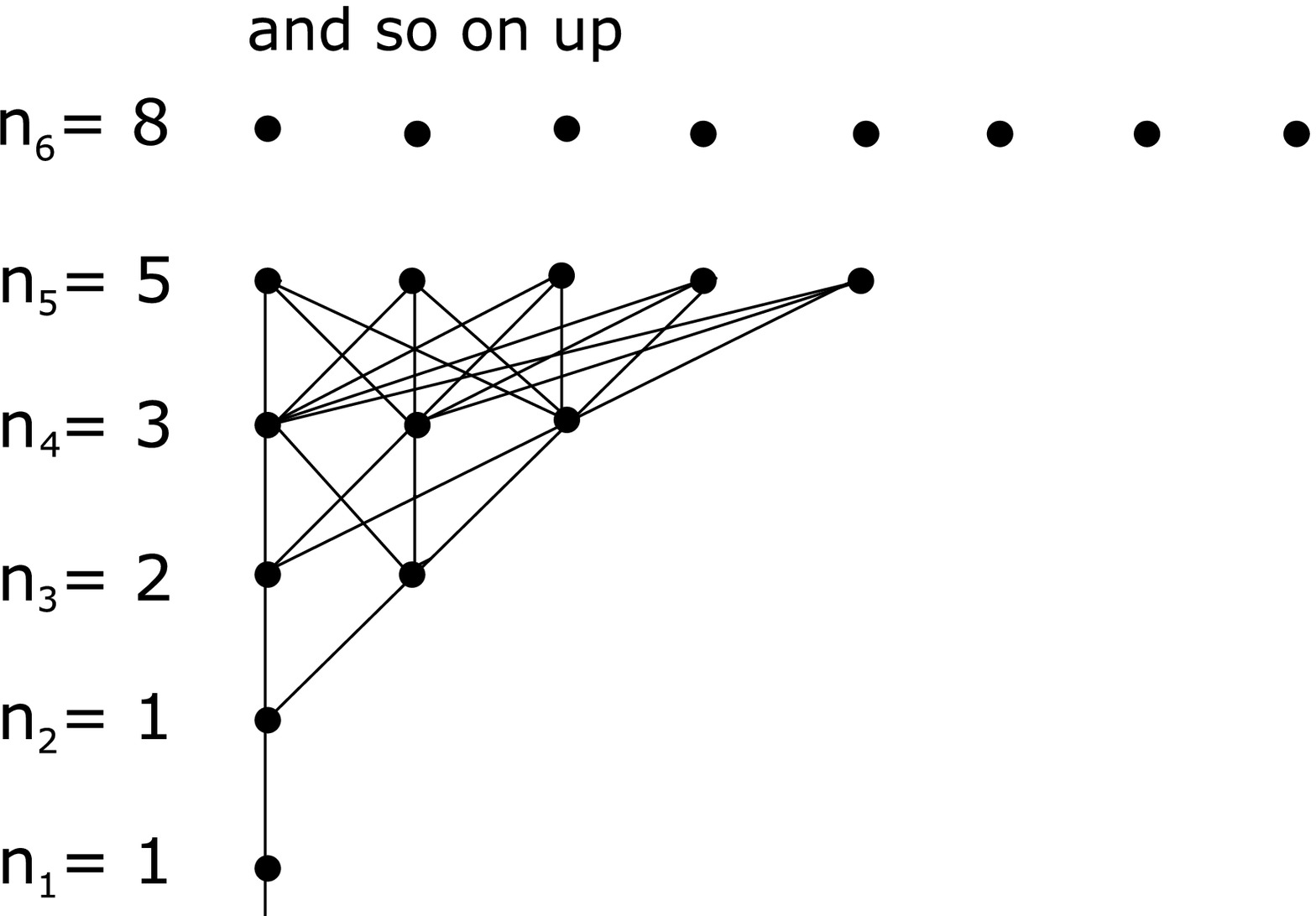}

\vspace{1mm}

\noindent {\small Fig.5. Display of Fibonacci numbers' cobweb poset.}
\end{center}

\noindent Compare this with  the bottom levels of a
Young-\textit{Fibonacci lattice}, introduced by Richard Stanley-
in Curtis Greene`s gallery of posets [2].

\vspace{2mm}

\noindent As seen above - for example  the $Fig.5$. displays the
rule of the construction of the  Fibonacci cobweb poset. It is
being visualized clearly while defining this \textit{non-lattice} cobweb 
poset  $P$ with help of its incidence matrix $\zeta$ [16]. The incidence 
matrix $\zeta$ representing uniquely just this cobweb poset $P$
exhibits (see below) a staircase structure characteristic to Hasse diagrams
of all cobweb posets. 

\vspace{1mm}

$$ \left[\begin{array}{ccccccccccccccccc}
1 & 1 & 1 & 1 & 1 & 1 & 1 & 1 & 1 & 1 & 1 & 1 & 1 & 1 & 1 & 1 & \cdots\\
0 & 1 & 1 & 1 & 1 & 1 & 1 & 1 & 1 & 1 & 1 & 1 & 1 & 1 & 1 & 1 & \cdots\\
0 & 0 & 1 & 1 & 1 & 1 & 1 & 1 & 1 & 1 & 1 & 1 & 1 & 1 & 1 & 1 & \cdots\\
0 & 0 & 0 & 1 & 0 & 1 & 1 & 1 & 1 & 1 & 1 & 1 & 1 & 1 & 1 & 1 & \cdots\\
0 & 0 & 0 & 0 & 1 & 1 & 1 & 1 & 1 & 1 & 1 & 1 & 1 & 1 & 1 & 1 & \cdots\\
0 & 0 & 0 & 0 & 0 & 1 & 0 & 0 & 1 & 1 & 1 & 1 & 1 & 1 & 1 & 1 & \cdots\\
0 & 0 & 0 & 0 & 0 & 0 & 1 & 0 & 1 & 1 & 1 & 1 & 1 & 1 & 1 & 1 & \cdots\\
0 & 0 & 0 & 0 & 0 & 0 & 0 & 1 & 1 & 1 & 1 & 1 & 1 & 1 & 1 & 1 & \cdots\\
0 & 0 & 0 & 0 & 0 & 0 & 0 & 0 & 1 & 0 & 0 & 0 & 0 & 1 & 1 & 1 & \cdots\\
0 & 0 & 0 & 0 & 0 & 0 & 0 & 0 & 0 & 1 & 0 & 0 & 0 & 1 & 1 & 1 & \cdots\\
0 & 0 & 0 & 0 & 0 & 0 & 0 & 0 & 0 & 0 & 1 & 0 & 0 & 0 & 1 & 1 & \cdots\\
0 & 0 & 0 & 0 & 0 & 0 & 0 & 0 & 0 & 0 & 0 & 1 & 0 & 1 & 1 & 1 & \cdots\\
0 & 0 & 0 & 0 & 0 & 0 & 0 & 0 & 0 & 0 & 0 & 0 & 1 & 1 & 1 & 1 & \cdots\\
0 & 0 & 0 & 0 & 0 & 0 & 0 & 0 & 0 & 0 & 0 & 0 & 0 & 1 & 0 & 0 & \cdots\\
0 & 0 & 0 & 0 & 0 & 0 & 0 & 0 & 0 & 0 & 0 & 0 & 0 & 0 & 1 & 0 & \cdots\\
0 & 0 & 0 & 0 & 0 & 0 & 0 & 0 & 0 & 0 & 0 & 0 & 0 & 0 & 0 & 1 & \cdots\\
. & . & . & . & . & . & . & . & . & . & . & . & . & . & . & . & . \cdots\\
 \end{array}\right]$$

\vspace{1mm} \noindent \textbf{Figure 6.  The incidence matrix
$\zeta$ for the Fibonacci cobweb poset}

\vspace{2mm}

\noindent \textbf{Note} that the knowledge of $\zeta$  matrix explicit
form enables one  to count via standard algorithms [16]
the M{\"{o}}bius matrix $\mu =\zeta^{-1} $ and other typical
elements of incidence algebra perfectly suitable for calculating
number of chains, of maximal chains etc. in finite sub-posets of
$P$. (All elements of the corresponding real incidence algebra are then
given by a matrix of the Fig.6 with $1$`s replaced by arbitrary real numbers).

\vspace{2mm}

\noindent \textbf{Cobweb posets as complete bipartite  digraph sequences.} 

\noindent As announced earlier the cobweb posets and their natural subposets $P_n$ are 
graded posets. They are vertex partitioned into  antichains $\Phi_k$ 
for $k = 0, 1, ..., r,...$ (where $r$ is a nonnegative integer) such that for each $\Phi_r$, all of the elements
covering $x$ are in $\Phi_{r+1}$ and all the elements covered by $x$ are in $\Phi_r$. We shall call
the $\Phi_n$ the  $n-th$-level.  $P_n$ is then $(n+1)$ level ranked poset. We are now in a position to 
observe that the cobweb posets may be identified with a chain of di-bicliques i.e. by definition - a chain of complete bipartite one direction digraphs. This is outstanding property as then any chain of relations is obtainable from the corresponding cobweb poset chain of complete relations just by deleting arcs in di-bicliques of this complete relations chain. Indeed.  For any given natural numbers valued sequence the graded (layered) cobweb posets` DAGs  are equivalently representations of a chain of binary relations. Every relation of the cobweb poset chain is bi-univocally represented by the uniquely designated  \textbf{complete} bipartite digraph - a digraph which is a di-biclique  designated by the very  given sequence. The cobweb poset may be therefore identified with a chain of di-bicliques i.e. by definition - a chain of complete bipartite one direction digraphs. Say it again, any chain of relations is obtainable from the cobweb poset chain of complete relations  via deleting  arcs (arrows) in chains` di-bicliques elements.\\ 

\noindent The above intuitively transparent names (a chain of di-cliques etc.) are the names
of the below  defined objects. These objects are natural correspondents of their undirected graphs relatives as
we can view a directed graph as an undirected graph with arrowheads added. Let us start with primary notions remembering that any cobweb subposet $P_k$ is a DAG of course. 

\noindent A bipartite digraph is a digraph whose vertices can be divided into two disjoint sets $V_1$ and $V_2$ such that every arc  connects a vertex in $V_1 $  and a vertex  in $V_2 $. Note that  there is no arc between two vertices in the same independent set  $V_1 $ or $V_2 $. No two nodes of the same partition set are adjacent.

\noindent A one direction bipartite digraph is a bipartite digraph such that every arc  originates at  a node in $V_1 $  and  terminates at a node in $V_2 $. The extension of "being one direction" to $k$-partite digraphs  is automatic. Note that  there is no arc  between two vertices in the same set.
Intuitively  - one may  color the nodes of a bipartite digraph black and blue such that no arc  exists between like colors.

\noindent A  $k$-partite not necessarily one direction  digraph $D$ is obtained from a $k$-partite undirected  $k$-partite 
graph $G$ by replacing every edge  ${xy}$  of  $G$  with the arc $\left\langle xy \right\rangle$, arc $\left\langle yx \right\rangle$  or both  $\left\langle xy \right\rangle$ and $\left\langle yx \right\rangle$. The partite sets of  $D$ are the partite sets of $G$.   
 
\begin{defn}. 
A simple directed graph $G = (V,E)$ is called bipartite if there exists a partition $V = V_1 + V_2$ of the vertex set
$V$  so that every edge (arc) in E is incident with  $v_1$  and  $v_2$ for some $v_1$  in $V_1$  and $v_2$ in $V_2$.
It is complete if  any node from $V_1$  is adjacent to all nodes of  $V_2$.
\end{defn}
We shall denote our special case biparte one direction digraphs as follows $G = (V_1+V_2,E) \equiv  L_{{F_k},F_{k+1}}$ to
inform that the partition has parts $V_1 $ and $V_2 $  where $\left|V_1 \right|= F_k$  and  $\left|V_2\right|= F_{k+1}.$

\noindent \textbf{Notation.} The  special $k$ - partite $\equiv$  $k$ - level one direction  digraphs considered here shall be
denoted  by the corresponding symbol $L_{p,q,...,r}$. The complete $k$ -partite $\equiv$  $k$ -level one direction complete digraph is coded  $K_{p,q,...,r}$  as in the non directed case (with no place for confusion because only one direction  digraphs are to be considered in what follows). 
\noindent $L_{{F_k},F_{k+1},...,F_n}\equiv \langle\Phi_k \rightarrow \Phi_n \rangle $ denotes  $(n-k+1)$ - partite one direction digraph  $\equiv (n-k+1)$ - level one direction digraph $ n\geq k$ whose partition has ( antichains from $\Pi$) the parts $\Phi_k, \Phi_{k+1},...,\Phi_n,\;  \left|\Phi_k \right|= F_k,  \left|\Phi_{k+1} \right|= F_{k+1},…,\left|\Phi_n \right|= F_n.$ Any  cobweb one-layer  $\langle\Phi_k \rightarrow \Phi_{k+1} \rangle,\quad k<n,\quad k,n \in N\cup\{0\}\equiv Z_\geq$  is a complete bipartite one direction digraph  $L_{{F_i},F_{i+1}} = K_{{F_i},F_{i+1}}$  i.e. by definition it is the \textit{di-biclique}. The cobweb one-layer vertices constitute bipartite set $V(\langle\Phi_k \rightarrow \Phi_{k+1} \rangle)=\Phi_k  + \Phi_{k+1}$ while edges are all those arcs incident with the two  antichains nodes in the poset $\Pi$ graph representation.\\
Any  cobweb subposet  $ P_n \equiv\langle\Phi_0 \rightarrow \Phi_n \rangle \equiv L_{F_0,F_1,...,F_n}, n\geq 0$  is a $(n+1)$-partite = $(n+1)$-level one direction digraph. We shall keep on calling  a complete bipartite one direction digraph a  \textbf{di-biclique} because it is a special kind of bipartite one direction digraph, whose every vertex of the first set is connected by an arc originated in this very node to \textbf{every} vertex of the second set of the given bi-partition. 
\noindent Any  cobweb layer $\langle\Phi_k \rightarrow \Phi_n \rangle,\quad k\leq n,\quad k,n \in N\cup\{0\}\equiv Z_\geq$  is a one direction  the $(n-k+1)$-partite $\equiv$ $(n-k+1)$-level one direction DAG  $L_{{F_k},F_{k+1},...,F_n}$ with an additional defining property: it is a chain of di-bicliques  for  $k<n$.  Since now on we shall identify both:  
                       $$\langle\Phi_k \rightarrow \Phi_n \rangle= L_{{F_k},F_{k+1},...,F_n}.$$
Note:  $P_n$ $\equiv$ $\langle\Phi_0 \rightarrow \Phi_n \rangle= L_{{F_0},F_{k+1},...,F_n},\; n\geq 0$  and 
$\langle\Phi_k \rightarrow \Phi_n \rangle \neg = K_{{F_k},F_{k+1},...,F_n}$   for $n>k+1.$\\
\noindent \textbf{Observation.}
$$ \left|E(\langle\Phi_k \rightarrow \Phi_{k+m} \rangle) \right|= {\sum}  ^{m-1}_{i=0}F_{k+i}F_{k+i+1},$$
where  $E(G)$ denotes the set of edges of a graph $G$ (arcs of a digraph $G$).

\vspace{2mm}

\noindent The following  property $(*)$

$$ (*)  \ \ \  \langle\Phi_k \rightarrow \Phi_{k+1}\rangle \equiv L_{{F_k},F_{k+1}}=K_{{F_k},F_{k+1}},\ \ k = 0,1,2,...$$
i.e. $L_{{F_k},F_{k+1}}$  is a di-biclique for $k = 0,1,2,...$ , might be considered the definition of  $F_0$ rooted 
$F$-cobweb graph $\Pi$ or in short $P$ if $F$-sequence has been established. The cobweb poset $\Pi$ is being thus identified with a \textbf{chain of di-bicliques}.  The usual convention is to  choose  $F_0 =1.$  One may relax this constrain, of course.

\noindent Thus any  cobweb one-layer  $\langle\Phi_k \rightarrow \Phi_{k+1}\rangle$ is a complete one direction  bipartite digraph    i.e.  a  di-biclique. This is how the definition of the $F$-cobweb graph $\Pi(F)$ or in short $P$ - has emerged.

\vspace{2mm}

\noindent For infinite cobweb poset $\Pi$  with the set of vertices $P =V(\Pi)$ one has an obvious $\delta(G)=2$ domatic vertex 
partition of this $F$-cobweb poset.  

\begin{defn}. 
A subset $D$ of the vertex set $V(G)$ of a digraph $G$ is called dominating in $G$, if each vertex of $G$ either is in $D$, or is adjacent to a vertex of $D$. Adjacent means that  there exists an originating or terminating arc in between the two- any node from $G$  outside $D$  and  a node from $D$.
\end{defn}

\begin{defn}. 
A domatic partition of $V$ is a partition of $V$ into dominating sets, and the number of these dominating sets is called the size of such a partition.  The domatic number  $ \delta(G)$ is the maximum size of a domatic partition.
\end{defn}

\vspace{2mm}

\noindent An infinite cobweb poset $\Pi$  with the set of vertices $P =V(\Pi)$  has the $\delta(\Pi)=2$ domatic vertex 
partition, namely a $mod \ 2$ - partition. 

\noindent $V = V_0 \cup V_1$ where  $V_0 = \bigcup _{k= 2s+1}\Phi_k , \  s = 0,1,2,... $  - ("black levels");
\noindent and  $V_1 = \bigcup _{k= 2s}\Phi_k , \  s = 0,1,2,... $  - ("blue levels").

\noindent Note:  Natural $mod \ n$ partitions of  the cobweb poset`s set of vertices $V(\Pi)=P$  ($n$ colours), 
$P = V_0 \cup V_1 \cup V_2 \cup ... \cup V_{n-1}$ ,  $V_i =\bigcup _{k= 2s+i} , s = 0,1,2,...,i \in Z_n = {0,1,...,n-1}$
for $n>2$ are \textbf{not domatic}.

\noindent  Cobweb one-layer or more than one-layer subposets  $\langle\Phi_k \rightarrow \Phi_{k+m}\rangle \equiv L_{{F_k},F_{k+1}}=K_{{F_k},F_{k+1}},\ \ k = 0,1,2,...$ have also correspondent, obvious the $ \delta(G)=2$ domatic partitions for $m>0.$

\vspace{4mm}

\noindent \textbf{Combinatorial interpretation.}

\vspace{2mm}

\noindent The crucial and elementary observation now is that an eventual
cobweb poset or any combinatorial interpretation of $F$-binomial
coefficients makes sense \textit{ not for arbitrary} $F$
sequences as $F-nomial$ coefficients should be nonnegative
integers (hybrid sets are not considered here).

\vspace{2mm}

\begin{defn}.
A natural numbers` valued sequence $F = \{n_F\}_{n\geq 0}$, $F_0 =1$ is called 
cobweb-admissible iff
$$ \left( \begin{array}{c} n\\k\end{array}\right)_{F}\in N_0\quad
for \quad k,n\in N_0.$$
\end{defn}
$F_0 = 0$ being sometimes acceptable as $0_F! \equiv F_0! = 1.$
\vspace{2mm}

\noindent Incidence coefficients of any reduced incidence algebra of full 
binomial type [16] immensely important for computer science are computed exactly with their correspondent cobweb-admissible sequences.
These include binomial (Newton) or  $q$- binomial  (Gauss) coefficients. For other $F$-nomial 
coefficients - computed with cobweb admissible sequences - see in what follows after Observation 3.

\vspace{2mm}

\noindent \textbf{Problem 1}. \textit{ Find effective characterizations  and/or
an algorithm to produce the cobweb admissible sequences i.e. find all examples.}

\vspace{2mm}

\noindent Right from the definition of $P$ via its
Hasse diagram pictures the important observations follow which lead to a specific, new joint
combinatorial interpretation of cobweb poset`s characteristic binomial-like coefficients.  

\vspace{3mm}
\noindent {\bf Observation 1}.

\noindent {\it The number of maximal chains starting from The Root  (level
$0_F$) to reach any point at the $n-th$ level  with $n_F$ vertices
is equal to $n_{F}!$}.

\vspace{2mm}

\noindent {\bf Observation 2}. $(k>0)$

\noindent {\it The number of all maximal chains in-between $(k+1)-th$ level $\Phi_{k+1}$
and the $n-th$ level $\Phi_n$ with $n_F$ vertices
is equal to $n_{F}^{\underline{m}}$, \quad where $m+k=n.$ } \\

\vspace{1mm}

\noindent Indeed. Denote the number of ways to get along maximal chains from 
\textit{any fixed point} (the leftist for example) in $\Phi_k$ to any vertex 
in  $\Phi_n , n>k$ with the symbol\\
  $$[\Phi_k \rightarrow \Phi_n]$$
then obviously we have ( $[\Phi_n \rightarrow \Phi_n]\equiv 1)$:\\
           $$[\Phi_0 \rightarrow \Phi_n]= n_F!$$ 
and
$$[\Phi_0 \rightarrow \Phi_k]\times [\Phi_k\rightarrow \Phi_n]=
[\Phi_0 \rightarrow \Phi_n].$$

\vspace{2mm}

\noindent For the purpose of a new joint combinatorial interpretation
of $F-sequence-nomial$ coefficients ({\it F-nomial} - in short)
let us consider all finite \textit{"max-disjoint"} sub-posets
rooted at the $k-th$ level at any fixed vertex $\langle
r,k \rangle, 1 \leq r \leq k_F $  and ending  at corresponding
number of vertices at the $n-th$ level ($n=k+m$) where the
\textit{max-disjoint} sub-posets are defined below.

\vspace{2mm}

\begin{defn}.
Two posets are said to be max-disjoint if considered as sets of maximal chains 
they are disjoint i.e. they have no maximal chain in common. An equipotent copy of $P_m$ 
[`\textbf{equip-copy}'] is  defined as such a maximal chains family  
equinumerous with $P_m$ set of maximal chains that the it constitutes a sub-poset
with one minimal element.
\end{defn}
We shall proceed with deliberate notation coincidence anticipating coming observation.
\begin{defn}.
Let us denote the number of all mutually max-disjoint equip-copies of $P_m$  rooted
at any  fixed vertex $\langle j,k \rangle , 1\leq j \leq k_F $ of $k-th$ level  with the symbol
$$ \left( \begin{array}{c} n\\k\end{array}\right)_{F}.$$
\end{defn}
One uses here the  customary convention:  $\left(
\begin{array}{c} 0\\0\end{array}\right)_{F}=1$ and $\left(
\begin{array}{c} n\\n\end{array}\right)_{F}=1.$

\vspace{2mm}
\noindent Compare the above with the Definition 4 and the Definition 10.

\vspace{1mm}

\noindent The number of ways to reach an upper level from a
lower one along any of  maximal chains  i.e.  the number of all
maximal chains from the level $\Phi_{k+1}$ to the level $\Phi_n ,\quad k>n$ 
is equal to 
  $$ [\Phi_k \rightarrow \Phi_n]= n_{F}^{\underline{m}}.$$

\noindent Therefore we have

\begin{equation}
\left( \begin{array}{c} n\\k\end{array}\right)_{F} \times [\Phi_0
\rightarrow \Phi_m] = [\Phi_k \rightarrow \Phi_n]=
n_{F}^{\underline{m}}
\end{equation}
where  $[\Phi_0 \rightarrow \Phi_m]= m_F!$ counts the number of
maximal chains in any equip-copy of  $P_m$. With this in mind we see
that the following holds.

\vspace{3mm}

\noindent {\bf Observation 3}. ${n,k}\geq 0$,

\noindent {\it Let $n = k+m$. Let $F$ be any cobweb admissible sequence.
Then the number of mutually max-disjoint  equip-copies i.e.  sub-posets
equipotent to $P_{m}$ , rooted at the same \textbf{ fixed} vertex of  $k-th$ level 
and ending at the n-th level is equal to}

$$\frac{n_{F}^{\underline{m}}}{m_{F}!} =
\left( \begin{array}{c} n\\m\end{array}\right)_{F}$$
$$ = \left(\begin{array}{c} n\\k\end{array} \right)_{F}=
\frac{n_{F}^{\underline{k}}}{k_{F}!}. $$

\vspace{2mm}

\noindent The immediate natural question now is

$$\Big\{{{\eta} \atop {\kappa}}\Big\}_{const}= ?$$ 
\vspace{2mm}
\noindent i.e.  the number of partitions with block sizes all equal to const = ?

\vspace{2mm}

\noindent where here $const=\lambda =  m_{F}!$ and 

$$ \eta =  n_{F}^{\underline{m}},\ \  \kappa = \left(\begin{array}{c} n\\k\end{array} \right)_{F} $$
\noindent The    \textit{const} indicates that this is the number of set partitions with  block sizes all equal to 
\textit{const} and we use Knuth notation $\Big\{{{\eta} \atop {\kappa}}\Big\}$ for  Stirling numbers of the second kind.

\vspace{2mm}

\noindent From the formula  (59) in [4] one infers  the Pascal-like matrix answer to the question above.

$$\Big\{{\eta \atop \kappa}\Big\}_{\lambda} = \delta_{\eta,\kappa \lambda}
\frac{\eta !}{\kappa !\lambda !^\kappa}.$$

\noindent  This gives us the rough upper bound for the number of tilings (see [1] for Pascal-like triangles) as we arrive now
to the following intriscically related problem.

\vspace{2mm}
\noindent \textbf{The partition or tiling Problem 2.}
\noindent Suppose now that  $F$  is a cobweb admissible sequence. Let us introduce  

$$\sigma P_m = C_m[F; \sigma <F_1, F_2,...,F_m>]$$
the equipotent sub-poset obtained from $P_m$ with help of a permutation $\sigma$
of the sequence $<F_1, F_2,...,F_m>$ encoding  $m$ layers
of $P_m$ thus obtaining the equinumerous sub-poset $\sigma P_m$ 
with the sequence $\sigma <F_1, F_2,...,F_m>$ encoding  now $m$ layers of $\sigma P_m$.
Then $P_m = C_m[F; <F_1, F_2,...,F_m>].$ Consider the layer 
$\langle\Phi_k \rightarrow \Phi_n \rangle,\quad k<n,\quad
k,n \in N$ partition into the equal size blocks which are here max-disjoint equi-copies of $P_m, m=n-k+1$.
The question then arises whether and under which conditions the layer may be
partitioned with help of max-disjoint  blocks of the form $\sigma P_m$. And how to visualize this phenomenon?
It seems to be the question of computer art, too. At first - we already know that an answer to 
the main question of such tilings existence - for some sequences $F$ -is in affirmative. 
Whether is it so for all cobweb admissible sequences -we do not know by now.
Some computer experiments done by student Maciej Dziemia\'nczuk [1] are encouraging. 
More than that. The author of [1]   proves tiling's existence for some cobweb-admissible sequences includin natural and Fibonacci numbers sequences. He show also that not all $F$ - designated cobweb posets do admit tiling as defined above. 
However problems:  "how many?" or "find it all" are opened. Let us recapitulate. \\

\vspace{2mm}

\noindent \textbf{Problem 2.  Recapitulation.}
Under which conditions layers may be partitioned with help of max-disjoint
blocks of established type  $\sigma P_m$ ? 
Find effective characterizations and/or find an algorithm to produce these partitions.

\noindent  The problem is intriguing also from the point of view of the art of progamming and smart  computer experimets.
One encounters such situation already looking for simple-minded rough upper bound for the number of tilings (see [1] for Pascal-like triangles).  The characterisctic feature and effect of calculations and experiments is inmediate appearence of Giant Numbers.

\vspace{2mm}

\noindent \textbf{GIGANTIC NUMBERS in COBWEB POSET TILING' upper bound PASCAL-LIKE TRIANGLES}

\noindent Examining in more detaile the answer to the already posed question above

$$\Big\{{{\eta} \atop {\kappa}}\Big\}_{const}= ?$$ 
we see that 

$$(*)\ \ \ \ \ \ \                       
\Big\{ {\eta \atop \kappa } \Big\}_\lambda = \delta_{\eta,\kappa\lambda}\Big\{ {\eta \atop \kappa } \Big\} = 
\delta_{\eta,\kappa\lambda} \frac{\eta!}{\kappa!}
\sum_{{{i_1 + i_2 + \ldots + i_\kappa = \eta} \atop {0 < i_1 = \ldots = i_\kappa = \lambda}}}^{}{\frac{1}{i_1! i_2!\ldots i_\kappa!}} = 
\delta_{\eta,\kappa\lambda}\frac{\eta!}{\kappa!(\lambda!)^\kappa}
$$
and we constatate  the simple-minded recurrence  (**) $\equiv$ (*) 
$$ (**)\ \ \ \ \ \ \ 
\Big\{ {{\eta + \lambda} \atop {\kappa + 1} } \Big\}_\lambda = \delta_{\eta+\lambda,(\kappa+1)\lambda}\frac{(\eta+1)^{\overline{\lambda}}}{\lambda!(\kappa+1)}
\Big\{ {\eta \atop \kappa} \Big\}_\lambda ;
\qquad 
\Big\{ {\lambda \atop 1} \Big\}_\lambda = 1
$$
resulting in

$$(***)\ \ \ \ \ \ \                       
\Big\{ {\kappa\lambda \atop \kappa } \Big\}_\lambda = \prod^{\kappa}_{k}\Phi_\lambda(k) \equiv [\Phi_\lambda(\kappa)]!$$
which we shall  call $\Phi_\lambda$ - factorial, where
$$ \Phi_\lambda(\kappa) = 
\frac{(\kappa\lambda)^{\underline{\lambda}}}{\lambda! \kappa}$$
\textbf{R\'esum\'e.}
The number of $\kappa$ - block partitions of $\kappa \lambda$ - numerous ensemble with block sizes all equal to
$\lambda =  m_{F}!$ is equal to its  $\Phi_\lambda$ - factorial.

\noindent \textbf{Experimental comparison}. Hereby  we show and compare few  upper bound numbers $\Phi_\lambda(\kappa)$ 
with the number $\Big\{ {n \atop k } \Big\}_F$  of \textbf{all} different tilings of the layer $\langle\Phi_k \rightarrow \Phi_n\rangle$ 
with  $\sigma P_m $blocks  as  computer experiments  show [1].
Except for obvious  cases  $\Phi_\lambda(\kappa)=1$  and very few first rows of upper-bound $F$ - triangles (zeros being not displayed), we are  facing   -  face to face -  enormously  Gigantic Numbers  as illustrated below via quotations from [1].

\begin{center}
\includegraphics[width=75mm]{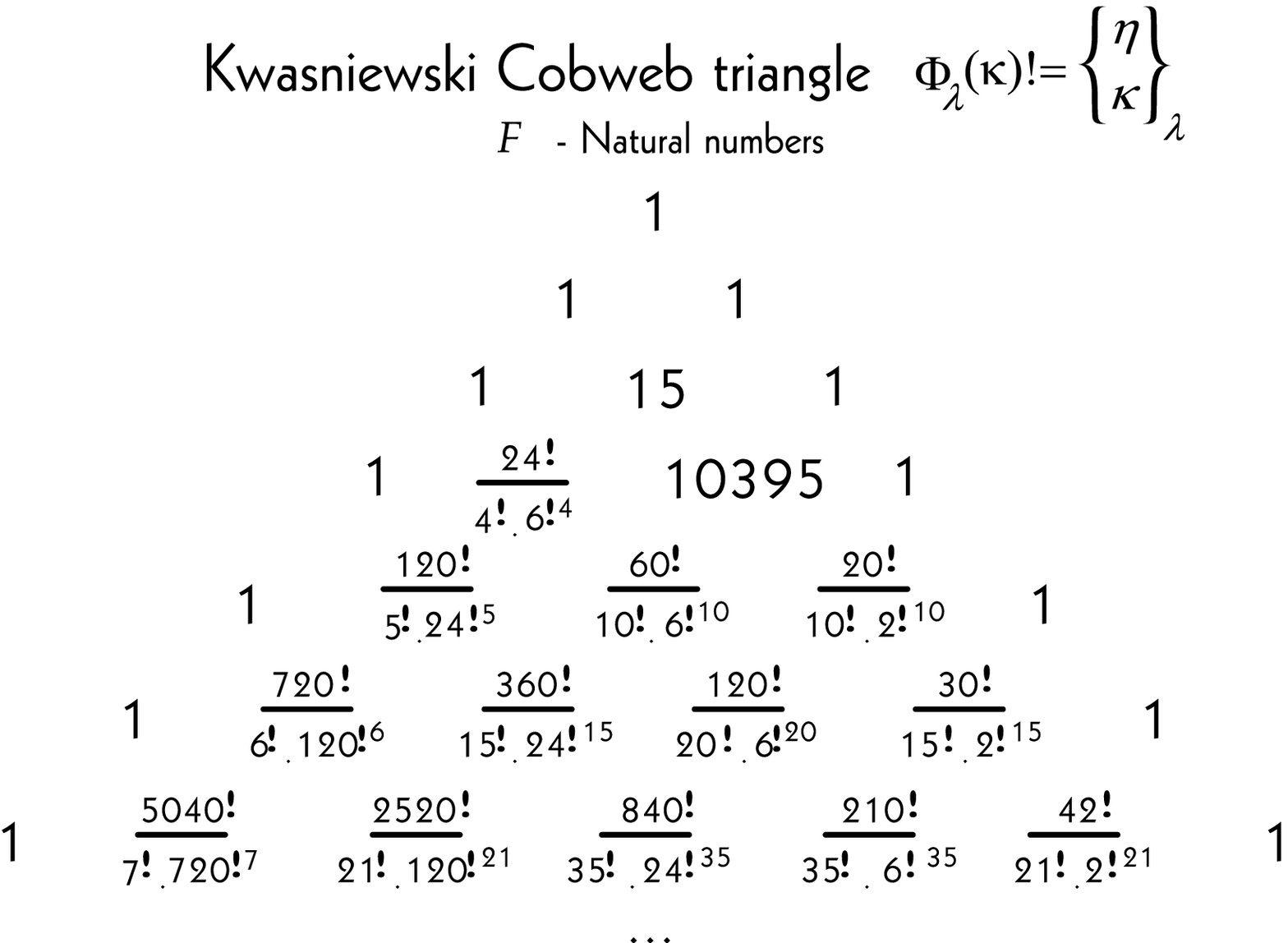}

\noindent {\small Fig.7. Display of Natural numbers' case - calculated upper bound.}
\end{center}

\begin{center}
\includegraphics[width=75mm]{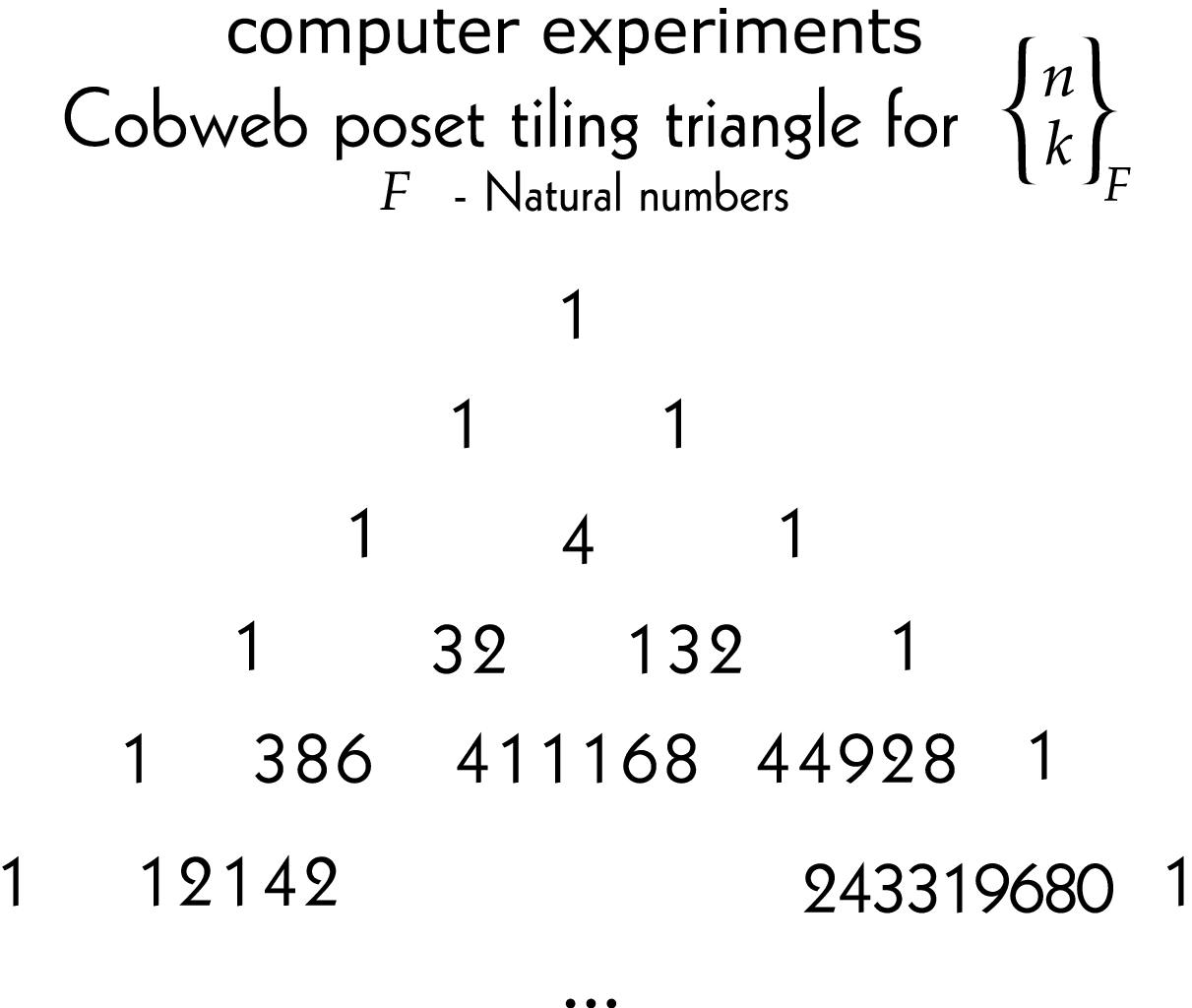}

\noindent {\small Fig.8. Display of Natural numbers case - number of tilings from experiment.}
\end{center}

\begin{center}
\includegraphics[width=75mm]{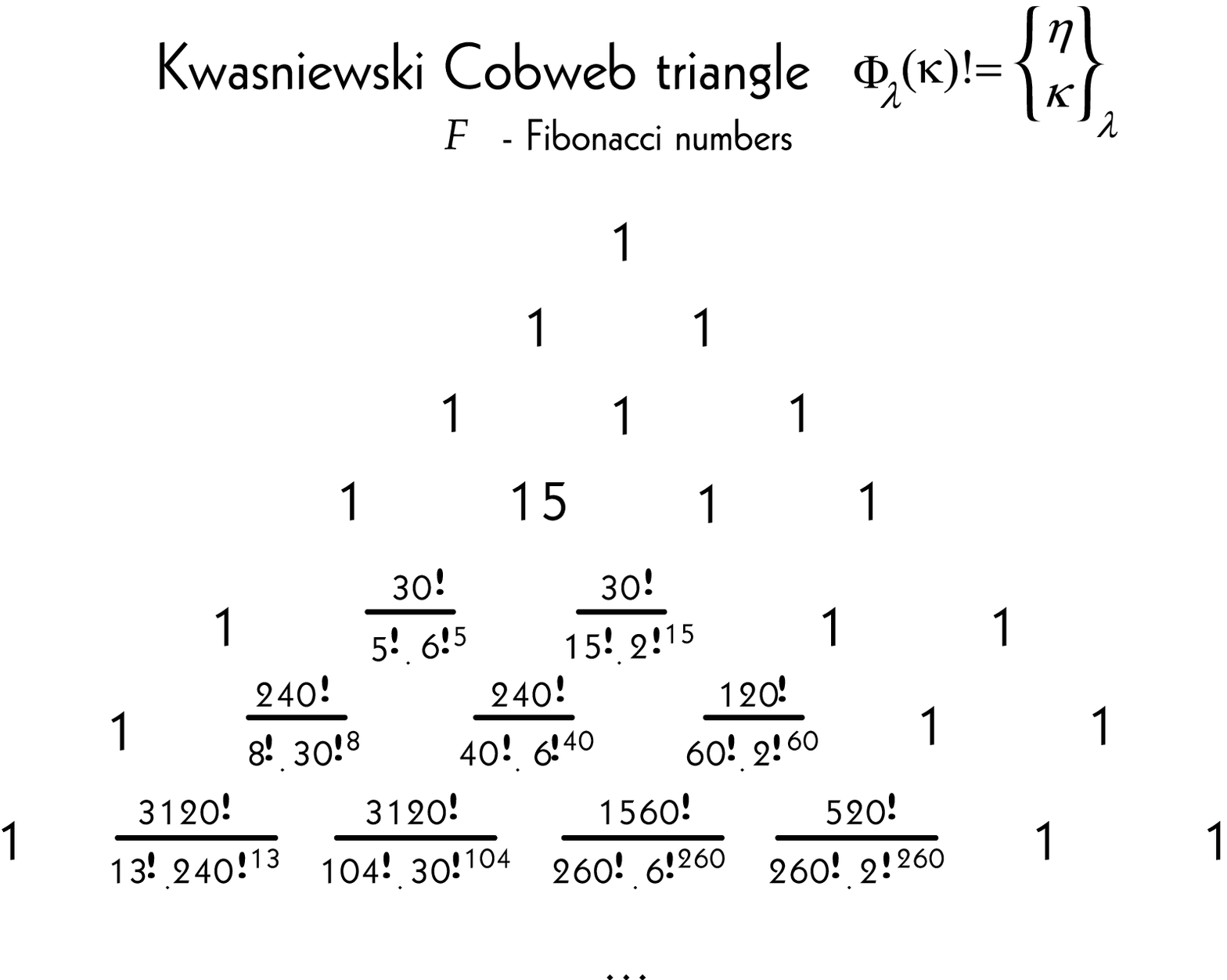}

\noindent {\small Fig.9. Display of Fibonacci numbers' case - calculated upper bound.}
\end{center}

\begin{center}
\includegraphics[width=75mm]{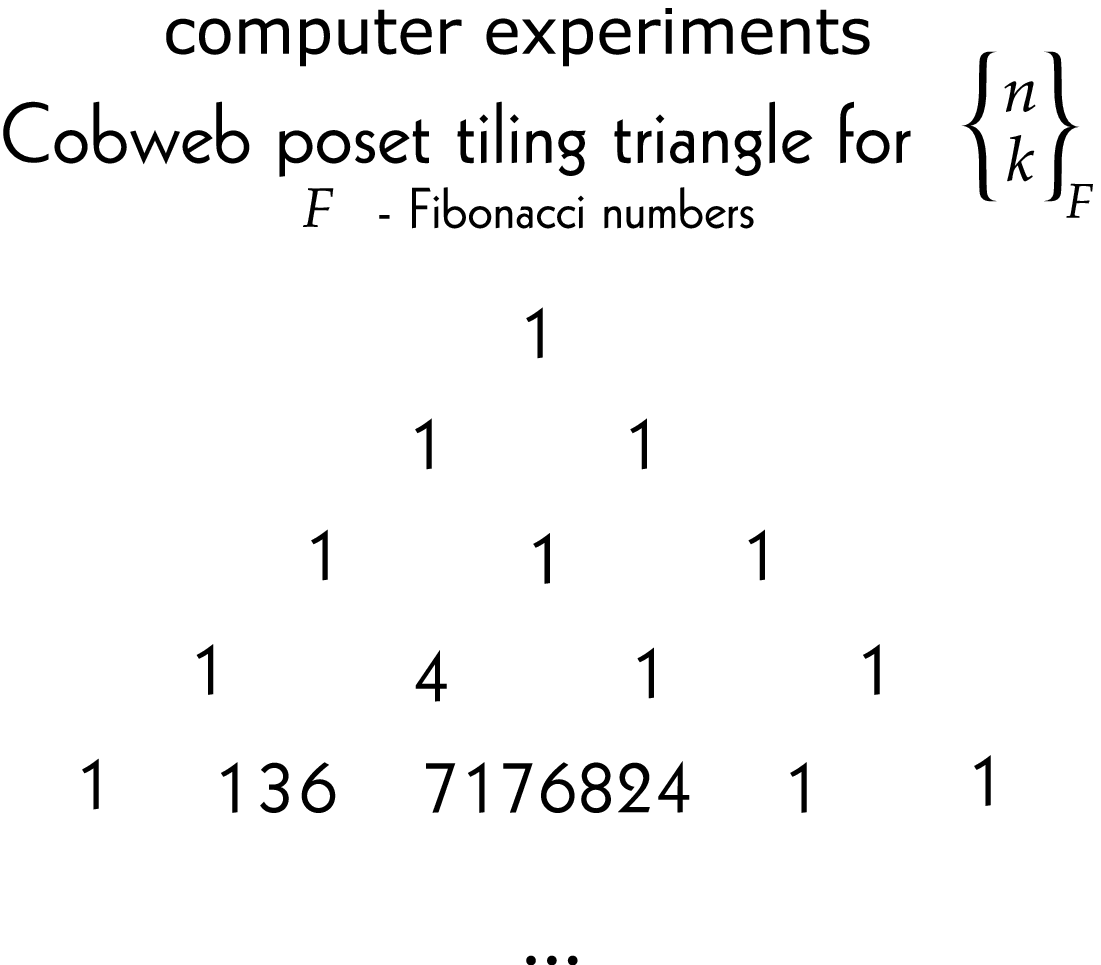}

\noindent {\small Fig.10.  Display of Natural numbers case - number of tilings from experiment.}
\end{center}

\vspace{3mm}

\noindent \textbf{Coming over to the last problem announced} let us note that the Observation 3. provides us
with the \textit{new} combinatorial interpretation of  the immense  class
of all classical $F-nomial$ coefficients including
binomial or Gauss $q$- binomial ones or Konvalina generalized binomial
coefficients of the first and of the second kind [4] - which include Stirling numbers too.
All these  $F$-nomial coefficients naturally are computed with their correspondent
cobweb-admissible sequences. More than that - the vast `umbral' family of  
$F$-sequences [10,11,12,13,6] includes also those which shall be called \textit{ "GCD-morphic"} sequences.
This means that $GCD[F_n,F_m] = F_{GCD[n,m]}$ where $GCD$ stays for Greatest Common Divisor. 
\begin{defn}.
The sequence of integers $F=\{n_F\}_{n\geq 0}$  is called the
GCD-morphic  sequence  if  $GCD[F_n,F_m] = F_{GCD[n,m]}$ where
$GCD$ stays for Greatest Common Divisor operator.
\end{defn}
The Fibonacci sequence is a much nontrivial [11,12,6] guiding example of GCD-morphic sequence.
Of course  \textit{not all }incidence coefficients of reduced incidence
algebra of full binomial type are computed with  GCD-morphic sequences
however these or that - if computed with the cobweb correspondent admissible sequences  
all are given the new, joint  cobweb poset combinatorial interpretation via Observation 3. 
More than that - in [8] a prefab-like combinatorial description of cobweb posets is being served with corresponding
generalization of the fundamental exponential formula. 

\noindent \textbf{Question:} which of these above mentioned sequences are GCD-morphic
sequences?

\vspace{2mm}

\noindent \textbf{GCD-morphism Problem. Problem III.} \textit{ Find effective characterizations
and/or  an algorithm to produce the GCD-morphic sequences i.e. find all examples.}

\vspace{2mm}

\noindent The recent papers on DAGs related to this article and its clue references [8,9] apart from [1] are [7] and [14].

\vspace{2mm}

\noindent \textbf{For the by now remark.}
\noindent For resent results  on cobweb posets  -up to the January 2008 - see [19]. The $poli-ary$  relations' ( as in relational database models)
point of view on cobweb posets and their KoDAGs  representatives    see  [20]  from 21 of December  2008  and for thus resulting natural open problems see then  [21].

\vspace{3mm}

\textbf{Acknowledgements} The author appraises much Maciej Dziemia\'nczuk's computer experiments
founded  on his vivid-active and effective interest in cobweb posets' now our joint investigations. I am 
also indebted for his successful indication of an important misprint. 
The author also appraises  Maciej Dziemia\'nczuk now Gda\'nsk Univerity Student's $TeX-nology$ aid.

The author expresses \textbf{solemnly}  his gratitude  to   Dr Ewa Krot-Sieniawska for her several years' cooperation
until \textbf{she was kicked out by Bialystok University authorities exactly on the day she had defended her Rota and KoDAs related dissertation with distinction}. Innocent was thus  penalized   \textbf{because of me} unmasking in public  misconducts,  infringements and contraventions.

\begin
{thebibliography}{99}
\parskip 0pt

\bibitem{1}
M. Dziemia\'nczuk,  {\it On cobweb posets tiling problem}, Adv. Stud. Contemp. Math. volume 16 (2), 2008 (April) pp. 219-233, 
$arXiv:math. Co/0709.4263, \ 4 \ Oct\   2007.$

\bibitem{2}
Curtis Greene, $www.haverford.edu/math/cgreene/posets/posetgallery.html$.  

\bibitem{3}
Hasse Diagram,  $http://mathworld.wolfram.com/HasseDiagram.html$. 

\bibitem{4}
Charles Jordan {\it On Stirling Numbers}  Tôhoku Math. J. \textbf{37} (1933),254-278.     

\bibitem{5}
J. Konvalina , {\it A Unified Interpretation of the Binomial
Coefficients, the Stirling Numbers and the Gaussian Coefficients}
The American Mathematical Monthly {\bf 107} (2000), 901-910.

\bibitem{6}
E. Krot, {\it An Introduction to Finite Fibonomial Calculus}, CEJM 2(5) (2005) 754-766.

\bibitem{7}
E. Krot, {\it The first ascent into the Fibonacci Cob-web Poset}, 
Adv. Stud. Contemp. Math.  \textbf{11} (2) (2007) 179-184.

\bibitem{8}
A. K. Kwa\'sniewski, {\it Cobweb posets as noncommutative prefabs}
Adv. Stud. Contemp. Math. \textbf{ 14 } (1) (2007) 37-47. 

\bibitem{9}
A. K. Kwa\'sniewski: {\em  First Observations on Prefab Posets' Whitney Numbers}
Advances in Applied Clifford Algebras\textbf{ Vol.18}, Number 1 / February, \textbf{2008} 57-73. arXiv: math.CO /0802.1696  cs.DM

\bibitem{10}
A. K. Kwa\'sniewski {\it Main  theorems of extended finite
operator calculus} Integral Transforms and Special Functions, {\bf
14} (6) (2003) 499-516.

\bibitem{11}
A. K. Kwa\'sniewski, {\it The logarithmic Fib-binomial formula},
Advanced Stud. Contemp. Math. {\bf 9} No 1 (2004) 19-26.

\bibitem{12}
A. K. Kwa\'sniewski, {\it Fibonomial cumulative connection constants},
Bulletin of the ICA  \textbf{44} (2005) 81-92.

\bibitem{13}
A. K.  Kwa\'sniewski, {\it On umbral extensions of Stirling numbers and Dobinski-like formulas}
Advanced Stud. Contemp. Math. {\bf 12}(2006) no. 1, pp.73-100.

\bibitem{14}
Anatoly D. Plotnikov, {\it About presentation of a digraph by dim 2 poset},
Adv. Stud. Contemp. Math.  \textbf{12} (1) (2006) 55-60

\bibitem{15}

Bruce E. SAGAN {\it Mobius Functions of Posets} (Lisbon lectures)IV: 
\textit{Why the Characteristic Polynomial factors }June 28 \textbf{2007} 
http://www.math.msu.edu/%7Esagan/Slides/mfp4.pdf

\bibitem{16}
E. Spiegel, Ch. J. O`Donnell  {\it Incidence algebras}  Marcel
Dekker, Inc., Basel, 1997.

\bibitem{17}

Richard P. STANLEY, \textit{Hyperplane Arrangements}, Proc. Nat.  Acad. Sci. \textbf{93} (1996), 2620-2625.
\textit{An Introduction to Hyperplane Arrangements } www.math.umn.edu/~ezra/PCMI2004/stanley.pdf

\bibitem{18}
M. Ward: {\em A calculus of sequences}, Amer.J.Math. \textbf{58} (1936)  255-266.

\bibitem{19}   A. Krzysztof Kwa\'sniewski, M. Dziemia\'nczuk, {\it   Cobweb posets - Recent Results},  Adv. Stud. Contemp. Math.  vol. 16 (2) April 2008 . pp. 197-218 ; arXiv: math. /0801.3985  Fri, 25 Jan 2008 17:01:28 GMT

\bibitem{20}  A. Krzysztof Kwa\'sniewski, {\it  Cobweb Posets and KoDAG Digraphs are Representing Natural Join of Relations, their diBigraphs and the Corresponding Adjacency Matrices}, arXiv:math/0812.4066v1 ,[v1] Sun, 21 Dec 2008 23:04:48 GMT 

\bibitem{21}  A. Krzysztof Kwa\'sniewski, {\it  Some Cobweb Posets Digraphs' Elementary Properties and Questions},  arXiv:0812.4319v1 ,[v1] Tue, 23 Dec 2008 00:40:41 GMT

\end{thebibliography}

%\begin{thebibliography}{99}
%{caban}P.~Caban, bleble

%\end{thebibliography}

\end{document}